\documentclass[12pt,a4paper,twoside]{article}
\usepackage[T2A]{fontenc}
\usepackage[utf8]{inputenc}
\usepackage[english]{babel}
\usepackage{amssymb, amsmath, amsfonts}
\usepackage[mathscr]{euscript}
\usepackage{mathrsfs}
\usepackage{theorem}
\usepackage{bm}
\usepackage[all]{xy}

\long\def\nodo#1{\relax}
\long\def\comment#1\endcomment{\relax}
\def\makepoint#1{\medbreak\noindent{\bf #1.\ }}
\def\zeropoint{\setcounter{subsection}{-1}}

\def\nxpoint{\refstepcounter{subsection}%
  \makepoint{\thesubsection}}
\def\nxsubpoint{\refstepcounter{subsubsection}%
  \makepoint{\thesubsubsection}}

\let\ptref=\refpoint

\let\myunderline=\mathbf
\let\stdfsets=\mathbf

\def\st#1{{\stdfsets{#1}}}
\def\stn{\st{n}}
\def\stm{\st{m}}
\def\stp{\st{p}}
\def\stq{\st{q}}

\def\card{\operatorname{card}}

\def\id{\operatorname{id}}

\def\Id{\operatorname{Id}}

\def\Ob{\operatorname{Ob}}

\def\Hom{\operatorname{Hom}}
\def\End{\operatorname{End}}

\def\Aut{\operatorname{Aut}}

\def\Ind{\operatorname{Ind}}
\def\Ker{\operatorname{Ker}}
\def\Coker{\operatorname{Coker}}

\def\Lim{\operatorname{Lim}}
\def\genLim#1{\operatorname{\vtop{\ialign{##\cr$\Lim$\cr\noalign{\nointerlineskip\kern0.5ex}\hfil$#1$\hfil\cr\noalign{\nointerlineskip\kern-0.5ex}\cr}}}}

\def\iHom{\operatorname{\myunderline{Hom}}}

\def\cU{\mathcal{U}}

\def\sA{\mathscr{A}}
\def\sB{\mathscr{B}}
\def\sC{\mathscr{C}}
\def\sE{\mathscr{E}}

\def\sG{\mathscr{G}}

\def\sQ{\mathscr{Q}}

\def\sS{\mathscr{S}}
\def\sT{\mathscr{T}}

\def\sZ{\mathscr{Z}}

\def\gS{\mathfrak{S}}

\def\simto{\stackrel\sim\to}
\def\Unit{{\bm 1}}
\def\bu{{\bm e}}

\let\phi=\varphi
\let\epsilon=\varepsilon

\let\emptyset=\varnothing
\let\injlim=\varinjlim
\let\projlim=\varprojlim

\let\textcat=\textsl

\def\catSets{{\textcat{Sets\/}}}

\def\catVectoid{{\textcat{Vectoid\/}}}

\def\catN{{\underline{\mathbb N}}}
\def\catAlg{\operatorname{\textcat{Alg\/}}}
\def\catCoalg{\operatorname{\textcat{Coalg\/}}}
\def\catCommAlg{\operatorname{\textcat{CommAlg\/}}}

\def\catFunct{\operatorname{\textcat{Funct\/}}}

\def\catEnd{\operatorname{\textcat{End\/}}}
\def\catHom{\operatorname{\textcat{Hom\/}}}
\def\catEndof{\operatorname{\textcat{Endof\/}}}

\def\catSheaves{\operatorname{\textcat{Sheaves\/}}}
\let\catDelta=\bbDelta

\def\catMod#1{{\textcat{$#1$-Mod}}}

\def\catVect#1{{\textcat{$#1$-Vect}}}

\theorembodyfont{\sl}
\newtheorem{Th}[subsection]{Theorem}

\newtheorem{Def}[subsection]{Definition}

\theorembodyfont{\rm}

\numberwithin{equation}{subsection}

\begin{document}
\title{Classifying Vectoids and Generalisations of Operads}
\author{Nikolai Durov\footnote{Supported by RFFR grant 08-01-00777-a}\\
\small{St.~Petersburg Department of Steklov Mathematical Institute}\\
\small{27 Fontanka emb., 191023 St.~Petersburg, Russia}}
\date{{}}

\maketitle

\begin{abstract}
A new generalisation of the notion of space, called {\em vectoid}, is
suggested in this work. Basic definitions, examples and properties are
presented, as well as a construction of direct product of vectoids.
Proofs of more complicated properties not used later are just
sketched.  Classifying vectoids of simplest algebraic structures, such
as objects, algebras and coalgebras, are studied in some detail
afterwards. Apart from giving interesting examples of vectoids not
coming from spaces known before (such as ringed topoi), monoids in the
endomorphism categories of these classifying vectoids turn out to
provide a systematic approach to construction of different versions of
the notion of an operad, as well as its generalisations, unknown
before.
\end{abstract}

\section{Definition of vectoid}

In this section we introduce basic definitions related to the notion of {\em vectoid\/}. Vectoid is a certain new generalisation of space, or rather ringed space (i.e.\ topological space with a sheaf of rings). In particular, any ordinary ringed space (e.g.\ an algebraic variety or a smooth manifold) determines a certain vectoid. At the same time vectoids generalise ordinary (not ringed) topological spaces and topoi as well.

The importance of vectoids reveals itself not only in the convenience of having such an all-encompassing theory, but is also due to the fact that many moduli spaces of algebraic structures (say, moduli space of associative algebras) don't exist or are extremely hard to construct in more conventional settings. At the same time such spaces (called {\em classifying vectoids}) can be rather easily constructed in the framework of vectoid theory.

Furthermore, endomorphisms of classifying vectoids for simple algebraic structures (such as algebras, coalgebras, pairs of algebras etc)\ turn out to be deeply related to operads and their modifications, thus providing a systematic approach to the study and construction of various operads and similar notions (called {\em algebrads} in this work), either known before or not.

The main idea of the way in which vectoids generalise ringed spaces is quite similar to the way in which topoi generalise topological spaces. Namely, in topos theory we associate to a topological space $X$ the category $\sE=\catSheaves(X)$ of sheaves of sets on~$X$, which is exactly the topos defined by this space, and then we try to work exclusively with this category $\sE$, and are even able to obtain an internal characterisation of such categories, not involving a topological space or site~$X$ (this is achieved by means of Giraud theorem, cf.\ \cite[IV.1]{SGA4}). Similarly, the vectoid corresponding to commutatively ringed space $(X,\sA)$ is just the category of sheaves of $\sA$-modules on~$X$, endowed with its natural tensor product $\otimes=\otimes_\sA$.

\begin{Def} {\em Prevectoid\/ $\sA$} is an ACU $\otimes$-category $\sA$, i.e.\ category $\sA$, endowed with associative and commutative tensor (or monoidal) structure $\otimes:\sA\times\sA\to\sA$, admitting a unit and satisfying the following conditions:
\begin{enumerate}
\item[V1)] Arbitrary (small) colimits exist in~$\sA$.
\item[V2)] Bifunctor $\otimes:\sA\times\sA\to\sA$ is {\em cocontinuous\/}, meaning that it preserves arbitrary colimits in any argument while the other argument is fixed: $X\otimes\injlim Y_\alpha\cong\injlim(X\otimes Y_\alpha)$.
\item[V3)] Finite limits exist in~$\sA$.
\item[V4)] Epimorphisms in $\sA$ are universally effective (i.e.\ are cokernels of their kernel pairs as well as their pullpacks). All equivalence relations in $\sA$ are efficient (i.e.\ any equivalence relation $R\subset X\times X$ coincides with $X\times_{X/R}X$).
\item[V5)] $\sA$ admits a small system of generators $\{G_\alpha\}$.
\end{enumerate}
\end{Def}

\nxsubpoint\label{p:prevect.term}
(Prevectoids and $\otimes$-prevectoids.)
In fact, it would be more natural to call this a {\em (commutative) $\otimes$-prevectoid}, and call any category satisfying V1) and V3)--V5) a prevectoid. We don't adopt such terminology here only because all (pre)vectoids we consider have a $\otimes$-structure.

\begin{Def} (Completion of a category.) Let $\sA$ be a category admitting arbitrary colimits. Denote by $\hat\sA:=\catFunct(\sA^{op},\catSets)$ the category of presheaves of sets on~$\sA$, and by $\tilde\sA\subset\hat\sA$ the full subcategory of~$\hat\sA$, consisting of {\em continuous\/} presheaves $F:\sA^{op}\to\catSets$, i.e.\ those which transform arbitrary colimits in $\sA$ into corresponding limits in $\catSets$:
\begin{equation}
F(\injlim_\alpha X_\alpha)\cong\projlim_\alpha F(X_\alpha)
\end{equation}
Identify each object $X\in\Ob\sA$ with corresponding representable functor $h_X:Y\mapsto\Hom_\sA(Y,X)$; by Yoneda $\sA$ is identified with a full subcategory of $\hat\sA$, and by definition of colimits any representable functor $h_X$ belongs to $\tilde\sA$. We obtain an embedding of full subcategories $\sA\subset\tilde\sA\subset\hat\sA$. We say that $\tilde\sA$ is the {\em completion\/} of category $\sA$.
\end{Def}

\begin{Def}
We say that category $\sA$ is {\em complete\/} (or {\em total\/} in an alternative terminology), if it admits arbitrary colimits, and any presheaf from $\tilde\sA$ is representable, i.e.\ the natural embedding $\sA\to\tilde\sA$ is an equivalence of categories.
\end{Def}

\begin{Def}\label{def:vectoid}
Prevectoid $\sA$ is called a {\em vectoid} if it satisfies the completeness axiom as well:
\begin{enumerate}
\item[V6)] $\sA$ is complete, i.e.\ $\sA\cong\tilde\sA$.
\end{enumerate}
\end{Def}

\nxpoint\label{p:vectoid.term} (Vectoids and $\otimes$-vectoids.)
A remark similar to \ptref{p:prevect.term} can be made here: it would be more natural to call the notion just defined a {\em (commutative) $\otimes$-vectoid}, and reserve the word `vectoid' for its variant without tensor structure, i.e.\ any category satisfying V1) and V3)--V6). However, we will continue to use terminology introduced above, and speak about {\em vectoids without tensor structure\/} when we don't insist on having a tensor structure.

\nxpoint (Comparison with topoi.)
The reader may have noticed that conditions V1)--V5) are similar to conditions of Giraud theorem (\cite[IV.1]{SGA4}), with the exception of universal disjointness of coproducts, which is generally false for vectoids. We introduce the completeness condition V6) instead. Roughly speaking, the counterpart of $\tilde\sA$ in topos theory is the category of sheaves with respect to canonical topology on some category $\sA$; thus the completeness condition essentially means that $\sA$ is something like a category of sheaves on itself. Recall that any topos $\sE$ is equivalent to the category $\tilde\sE$ of sheaves over itself with respect to the canonical topology; cf.\ \cite[IV]{SGA4}).

\nxpoint (Existence of arbitrary limits in vectoids.)
It is easy to see that V1) and V6), i.e.\ completeness of~$\sA$, imply existence of arbitrary (small) limits. Indeed, any limit of continuous presheaves is continuous, i.e.\ $\tilde\sA\subset\hat\sA$ is stable under arbitrary limits. Since V6) means that all continuous presheaves are representable, we see that arbitrary limits exist in $\sA\cong\tilde\sA$. In particular, condition V3) turns out to be superfluous; in fact, it is needed only to be able to state V4) for arbitrary prevectoids, not necessarily complete.

\nxpoint\label{p:clo.tens.vect} (Closedness of vectoid tensor structure.)
The {\em closedness} of tensor structure $\otimes$ of any vectoid~$\sA$, i.e. existence of inner $\iHom:\sA^{op}\times\sA\to\sA$, characterised by existence of canonical isomorphism $\Hom(X,\iHom_\sA(Y,Z))\cong\Hom(X\otimes Y,Z)$, can be verified in a similar fashion. Indeed, fix arbitrary objects $Y$, $Z\in\Ob\sA$; then functor $X\mapsto\Hom(X\otimes Y,Z)$ is continuous by V2), hence representable by some object $\iHom(Y,Z)$ according to~V6).

\begin{Def} (Vectoid morphisms.)
Let $\sA=(\sA,\otimes_\sA)$ and $\sB=(\sB,\otimes_\sB)$ be two vectoids.
A {\em morphism\/} $f$ of vectoid $\sA$ into vectoid $\sB$ is by definition a couple of adjoint
functors $(f^*,f_*)$, $f^*:\sB\to\sA$, $f_*:\sA\to\sB$, where the ``inverse image functor'' $f^*:\sB\to\sA$ is a $\otimes$-functor, i.e.\ is equipped with functorial isomorphisms $f^*(X\otimes_\sB Y)\cong f^*X\otimes_\sA f^*Y$, 
compatible with associativity, commutativity and unity constraints in the standard fashion.
\end{Def}

\nxpoint (2-category of vectoids.)
Similarly to the topoi, the collection of all vectoids and morphisms between them constitutes not just a category, but a 2-category. In other words, if $\sA$ and $\sB$ are vectoids, then $\Hom(\sA,\sB)$ is not just a set, but a category. A morphism between two objects $f$, $g$ of this category is defined to be a natural transformation of ``inverse image'' $\otimes$-functors $\theta:f^*\to g^*$. Notice that it induces a natural transformation of ``direct image functors'' $\theta':g_*\to f_*$ in the opposite direction. Nevertheless, we prefer to define 2-morphisms (i.e.\ morphisms between morphisms) in terms of inverse image functors, because it turns out to be more convenient while studying classifying vectoids.

\nxpoint (Vectoid morphism is determined by the inverse image functor.)
Notice that the direct image functor $f_*:\sA\to\sB$ is determined by the inverse image functor $f^*:\sB\to\sA$ uniquely up to an isomorphism, beeing its right adjoint. The inverse image functor $f^*$ is a $\otimes$-functor (by definition of vectoid morphism) and cocontinuous (commuting with arbitrary colimits), since it admits a right adjoint $f_*$. 

Conversely, any cocontinuous $\otimes$-functor $f^*:\sB\to\sA$ admits a right adjoint $f_*$ (by a reasoning similar to~\ptref{p:clo.tens.vect}), hence defines a vectoid morphism $f=(f^*,f_*):\sA\to\sB$, unique up to an isomorphism. Therefore, it suffices to construct cocontinuous $\otimes$-functor $f^*$ to determine vectoid morhism~$f$.

\section{Examples of vectoids}

\zeropoint\nxpoint
In this section we consider only those examples of vectoids which are constructed from already known spaces or topoi. Examples of vectoids which cannot be obtained in such a fashion will be presented later.

\nxpoint (Topological spaces.)
Any topological space $X$ determines a vectoid $(\catSheaves(X),\times)$, which is simply the category of sheaves of sets on~$X$ (i.e.\ the topos corresponding to~$X$) with $\otimes$-structure given by direct product of sheaves of sets ($\otimes=\times$). Continuous map $f:X\to Y$ yields a corresponding vectoid morphism $f=(f^*,f_*):\catSheaves(X)\to\catSheaves(Y)$ in the same direction. This is the reason why vectoid morphisms were defined so as to have the same direction as the direct image functors $f_*$, not the inverse image functors $f^*$.

\nxpoint (Topoi.)
More generally, any topos $\sE$ is a vectoid, if we consider the $\otimes$-structure given by direct products of objects ($\otimes=\times$). Thus we obtain a 2-functor which ``embeds'' 2-category of topoi into 2-category of vectoids. Notice that this functor needn't be fully faithful (i.e.\ inducing equivalence on each category of morphisms): a vectoid morphism $f:\sE\to\sE'$ between topoi $\sE$ and $\sE'$ is not necessarily a topos morphism since in topos theory $f^*$ is required to be cocontinuous and left exact (i.e.\ commuting with finite limits), while in the case of vectoid morphisms $f^*$ must be cocontinuous and preserve direct products, which is a somewhat weaker condition.

\nxpoint (Sheaf and presheaf categories.)
In particular, if $\sS$ is a Grothendieck site, then the sheaf category $\tilde\sS$ is a topos, hence a vectoid as well. If $\sS$ is an arbitrary small category, then the presheaf category $\hat\sS=\catFunct(\sS^{op},\catSets)$ is a topos and a vectoid. That's why each time we consider a presheaf category $\hat\sS$ without explicitly mentioning the choice of tensor product $\otimes$, we always mean $\otimes=\times$.

\nxpoint (Simplicial and cosimplicial sets.)
Taking point category, category of non-empty finite ordinals $\catDelta$ or its opposite $\catDelta^{op}$ as $\sS$, we obtain 
the category of sets $\catSets$, simplicial sets $s\catSets$ and cosimplicial sets $c\catSets$ as the corresponding category of presheaves.

\nxpoint (Final vectoid.) 
Notice that $(\catSets,\times)$ is 2-final object in the 2-category of vectoids, i.e.\ the category $\Hom(\sA,\catSets)$ is equivalent to the point category for any vectoid~$\sA$. 

\nxpoint (Ringed spaces.)
Now let $(X,\sA)$ be a commutatively ringed topological space (almost all kinds of varieties and manifolds studied in geometry have such a structure). Then the category of sheaves of $\sA$-modules $({\catMod\sA}_X,\otimes_\sA)$, endowed with its natural tensor product, is a vectoid. Morphisms of ringed spaces induce morphisms of corresponding vectoids.

\nxpoint (Ringed topoi.)
The previous example is immediately generalised to the case of commutatively ringed topos $(\sE,\sA)$: in this case $({\catMod\sA}_\sE, \otimes_\sA)$ (cf.\ \cite[IV]{SGA4}) is a vectoid.

\nxpoint (Rings.)
We can also consider the case of a sheaf of rings over a point, i.e.\ ringed space $(*,A)$, where $A$ is now an ordinary commutative ring. In this case we obtain the category of $A$-modules $\catMod{A}$ with the usual tensor product; thus $(\catMod{A},\otimes_A)$ is a vectoid. Notice that the 2-functor from commutative rings into vectoids just constructed is contravariant.

\nxpoint (Generalised rings.)
If $\Sigma$ is a {\em generalised ring\/} in the sense of~\cite{Du} (i.e.\ a commutative algebraic monad), then $(\catMod\Sigma,\otimes_\Sigma)$ is a vectoid. Thus the previous example can be extended to generalised rings.

\nxpoint\label{p:gen.rg.top} (Generalised ringed spaces and topoi.)
Similarly, categories of (sheaves of) modules over generalised ringed spaces and topoi are vectoids.

\nxpoint (Unit-generated vectoids.)
Let $\sA$ be a vectoid generated by its unit object $\Unit_\sA$. Suppose that this object is finitely presented (i.e.\ $\Hom_\sA(\Unit_\sA,-)$ commutes with filtered inductive limits) and projective (i.e.\ $\Hom_\sA(\Unit_\sA,-)$ transforms strict epimorphisms into surjective maps). Then $\sA=(\catMod\Sigma,\otimes_\Sigma)$ for a certain uniquely determined generalised ring $\Sigma$. Namely, let $p_*:=\Hom_\sA(\Unit_\sA,-):\sA\to\catSets$, $p^*:S\mapsto \Unit_\sA^{(S)}$ be its left adjoint; then $\Sigma:=p_*p^*$ is a monad on $\catSets$, algebraic because $\Unit_\sA$ is finitely presented. General results on algebraic theories due to Lawvere (см.\ \cite{Lawvere}) imply that the comparison functor $\sA\to\catSets^\Sigma=\catMod\Sigma$ is an equivalence of categories (here condition V4) turns out to be useful), and existence of tensor product on $\catMod\Sigma\cong\sA$, such that $\Sigma(1)\otimes\Sigma(1)\cong\Sigma(1)$, implies commutativity of~$\Sigma$. 

\nxpoint (Proof of completeness of $\catMod\Sigma$.)
Validity of conditions V1)--V5) in all examples listed above is obvious or at least well-known (the case of generalised rings can be found in \cite[4]{Du} and \cite[5]{Du}). However, the completeness condition V6) is less trivial. For the topos case it suffices to observe that cocontinuous presheaves of sets on topos $\sE$ are sheaves with respect to the canonical topology (since they commute with cokernels of kernel pairs of strict epimorphisms), hence are representable according to \cite[IV]{SGA4}. The commutative ring case might be treated directly; however, since any commutative ring is a generalised ring as well, it suffices to treat the generalised ring case. Combining the proofs for the generalised ring case and the topos case, one can obtain long and quite technical, but not illuminating proof of completeness for the case \ptref{p:gen.rg.top}, which encompasses all our other examples.

Therefore, we will content ourselves with the completeness proof for the case of $\sA=(\catMod\Sigma,\otimes_\Sigma)$, where $\Sigma$ is an arbitrary generalised ring.

\nxsubpoint (Notations related to generalised rings.)
So, fix a generalised ring $\Sigma$, i.e.\ a commutative algebraic monad $\Sigma=(\Sigma,\mu,\epsilon)$ on the category of sets $\catSets$. This means that $\Sigma:\catSets\to\catSets$ is an endofunctor commuting with filtered inductive limits, $\mu:\Sigma^2\to\Sigma$ and $\epsilon:\Id_\catSets\to\Sigma$ are some natural transformation subject to certain conditions (associativity and unit relations). Let $\sA:=\catMod\Sigma=\catSets^\Sigma$; recall (cf.\ \cite[4]{Du}) that $\catMod\Sigma$ consists of $\Sigma$-modules $M=(M,\alpha)$, i.e.\ sets $M$ together with a $\Sigma$-action $\alpha:\Sigma(M)\to M$, subject to some conditions. Let $\Gamma:\sA\to\catSets$ be the forgetful functor $(M,\alpha)\to M$, $L:\catSets\to\sA$ be its left adjoint, given by $L:S\mapsto(\Sigma(S),\mu_S)$; $\Sigma$-module $L(S)$ is said to be the {\em free $\Sigma$-module generated by set~$S$}, and is usually denoted by $\Sigma(S)$. For any integer $n\geq0$ we denote by $\stn:=\{1,2,\ldots,n\}$ the standard $n$-element set, and we put $\Sigma(n):=\Sigma(\stn)$; this set with the natural $\Sigma$-action $\mu_\stn$ is called the {\em free $\Sigma$-module with $n$ generators}.

Recall that $\Hom_\sA(\Sigma(S),M)\cong\Hom_\catSets(S,\Gamma(M))=\Gamma(M)^S=M^S$; in particular, $\Gamma(M)=\Hom_\sA(\Sigma(1),M)=\Hom_\sA(\Unit_\sA,M)$. Tensor product on $\sA$ is defined so as to have $\Sigma(S)\otimes\Sigma(T)\cong\Sigma(S\times T)$ and it is required to be right exact; it is uniquely determined by these two conditions since any $\Sigma$-module can be written as a cokernel of a couple of arrows between free modules.

\nxsubpoint (Properties of a continuous presheaf on~$\sA$.)
Now let $F:\sA^{op}\to\catSets$ be any continuous presheaf of sets on~$\sA$.
If it is representable by some $\Sigma$-module $M=(M,\alpha)$, then
$M=\Hom_\sA(\Unit_\sA,M)=F(\Unit_\sA)$ as a set; hence it makes sense to begin our proof of representability of $F$ by putting $M:=F(\Unit_\sA)$, and continue by constructing a $\Sigma$-action $\alpha:\Sigma(M)\to M$.

Denote $\Gamma(U):=U(\Unit_\sA)$ for any presheaf of sets $U$ on~$\sA$. Resulting functor $\Gamma:\hat\sA\to\catSets$ extends previously defined functor $\Gamma:\sA\to\catSets$. We have $M=\Gamma(F)$ for our chosen continuous presheaf of sets~$F$.

\nxsubpoint ($\Sigma$-action on $M$.)
Let's construct the action $\alpha:\Sigma(M)\to M$. Fix any $t\in\Sigma(M)$ for this. Since $\Sigma(M)\cong\Hom_\sA(\Unit_\sA,\Sigma(M))$, $t$ defines a morphism $\tilde t:\Unit_\sA\to L(M)=\Sigma(M)$ in~$\sA$. Moreover, functor $L$ is cocontinuous, hence $\Sigma(M)=\Unit_\sA^{(M)}$ (infinite coproduct in~$\sA$). Functor $F$ is continuous, hence $F(\Sigma(M))=F(\Unit_\sA^{(M)})=F(\Unit_\sA)^M=M^M$; thus $F(\tilde t):F(\Sigma(M))\to F(\Unit_\sA)$ can be considered as a map of set $M^M$ into $M$. However, the set $M^M=\Hom_\catSets(M,M)$ contains a distinguished element $\id_M$; put
\begin{equation}
\alpha(t):=(F(\tilde t))(\id_M)\in M
\end{equation}
We have thus constructed a map $\alpha:\Sigma(M)\to M$.

\nxsubpoint (Generalisation of previous construction.)
In fact, this construction works starting from any $t\in\Sigma(S)$, where $S$ is an arbitrary set now, yielding a map $F(\tilde t):M^S\to M$. This construction is functorial in~$S$, i.e.\ for any $\phi:S\to T$, $s\in\Sigma(S)$ and $z\in M^T=\Hom_\catSets(T,M)$ we have
\begin{equation}
F(\tilde s)(z\circ\phi)=F(\tilde t)(z),\quad\text{where $t=(\Sigma(\phi))(s)$}
\end{equation}

\nxsubpoint (Compatibility of $\alpha$ with the monad unit $\epsilon$.)
Applying this to $S=\st1=\{1\}$, $T=M$, $s=\bu=\epsilon_1(1)\in\Sigma(1)$, $\phi=i_x:1\mapsto x\in M$, $z=\id_M$, we see that $t=(\Sigma(i_x))(\epsilon_1(1))=\epsilon_M(x)$, since $\epsilon:\Id_\catSets\to\Sigma$ is a natural transformation, hence $\alpha(\epsilon_M(x))=\alpha(t)=F(\tilde t)(\id_M)=F(\tilde s)(i_x)=F(\id_{\Sigma(1)})(i_x)=i_x=x$ if we identify $M^1$ with $M$. Since $x\in M$ was arbitrary, this means that $\alpha\circ\epsilon_M=\id_M$, i.e.\ the map $\alpha$ is compatible with the monad unit~$\epsilon$.

\nxsubpoint (Compatibility of $\alpha$ with the monad multiplication $\mu$.)
Let's show now that $\alpha\circ\Sigma(\alpha)=\alpha\circ\mu_M:\Sigma^2(M)\to M$; together with our previous results this will mean that $(M,\alpha)$ is a $\Sigma$-module. Fix arbitrary $s\in\Sigma^2(M)$ for this, let $t:=\mu_M(s)\in\Sigma(M)$, and consider the following diagram in~$\hat\sA$:
\begin{equation}\label{diag:mult.comp}
\xymatrix@+1pc{
\Unit_\sA\ar[r]_{\tilde s}\ar@/^2pc/[rr]^{\tilde t}&
\Sigma^2(M)\ar[r]_{\mu_M}\ar[rd]^<>(.7){\alpha^+}&
\Sigma(M)\ar[d]^{\id_M^+}\\
&\Unit_\sA\ar[u]^{\Sigma(i_u)}\ar[ur]^<>(.3){\tilde u}&F}
\end{equation}
Here $\tilde x:\Unit_\sA=\Sigma(1)\to\Sigma(X)$ denotes the morphism in $\sA$ corresponding to some element $x\in\Sigma(X)$, for arbitrary set $X$; $i_x:\st1=\{1\}\to X$ denotes the embedding of one-element set into set $X$ with image $x$; finally, $\beta^+:\Sigma(X)\to F$ denotes the element of $\Hom_{\hat\sA}(\Sigma(X),F)=F(\Sigma(X))=F(\Unit_\sA^{(X)})=F(\Unit_\sA)^X=\Hom_\catSets(X,M)$, corresponding to a map of sets~$\beta:X\to M$. Apart from that, an arbitrary element $u\in\Sigma(M)$ is used in this diagram.

Notice that $\tilde t=\mu_M\circ\tilde s$ since $t=\mu_M(s)$;
equality $\mu_M\circ\Sigma(i_u)=\tilde u$ is also evident (both sides map the only generator $\bu=\epsilon_1(1)\in\Sigma(1)$ into $u$). In order to check $\id_M^+\circ\mu_M=\alpha^+$ inside $\Hom_{\hat\sA}(\Sigma^2(M),F)\cong\Hom_{\catSets}(\Sigma(M),M)$ it suffices to verify that the image of $\id_M^+\circ\mu_M$ in $\Hom_{\catSets}(\Sigma(M),M)$ coincides with $\alpha$, i.e.\ that its value on arbitrary $u\in\Sigma(M)$ equals $\alpha(u)$. Indeed, this value is $(\id_M^+\circ\mu_M)\circ\Sigma(i_u)=\id_M^+\circ\tilde u=\alpha(u)$ by definition of~$\alpha$. This completes the verification of commutativity of the above diagram.

Now notice that $\alpha^+=\id_M^+\circ\Sigma(\alpha)$ (this is a special case of general equality $(\beta\circ\phi)^+=\beta^+\circ\Sigma(\phi)$, valid for any $\beta:X\to M$ and $\phi:X'\to X$), whence $\alpha^+\circ\tilde s=\id_M^+\circ\Sigma(\alpha)\circ\tilde s=\id_M^+\circ\widetilde{\Sigma(\alpha)(s)}=\alpha\bigl(\Sigma(\alpha)(s)\bigr)=(\alpha\circ\Sigma(\alpha))(s)$.
On the other hand, commutativity of the diagram implies $\alpha^+\circ\tilde s=\id_M^+\circ\mu_M\circ\tilde s=\id_M^+\circ\tilde t=\alpha(t)=\alpha(\mu_M(s))$; this completes the proof of $\alpha\circ\Sigma(\alpha)=\alpha\circ\mu_M$.

\nxsubpoint (Construction of morphism $\theta:M\to F$.)
So we have constructed an object $M=(M,\alpha)$ of category $\sA$; let us construct a natural morphism $\theta:M\to F$ as well. By Yoneda this amounts to constructing an element $\theta\in F(M)$. Since for any object of $\catMod\Sigma$, for $M$ in particular, the following diagram is right exact:
$$\xymatrix{\Sigma^2(M)\ar[rr]<2pt>^{\mu_M}\ar[rr]<-3pt>_{\Sigma(\alpha)}&&
\Sigma(M)\ar[rr]^{\alpha}&&M}\quad,$$
we obtain $F(M)=\Ker\bigl(F(\Sigma(M))\rightrightarrows F(\Sigma^2(M))\bigr)=\Ker(M^M\rightrightarrows M^{\Sigma(M)})$, functor $F$ being continuous. Now notice that $\id_M\in M^M$ belongs to this kernel: this was actually shown during the verification of commutativity of diagram~\eqref{diag:mult.comp}. Hence we can take this element as~$\theta$.

\nxsubpoint ($\theta$ is isomorphism.)
It remains to show that $\theta:M\to F$ is a presheaf isomorphism, i.e.\ that for any $X\in\Ob\sA$ the map $\theta_X:\Hom_\sA(X,M)\to F(X)$ is bijective. Since any object $X$ may be written as cokernel of a couple of arrows between free $\Sigma$-modules, and $F$ is continuous hence transforms cokernels into kernels, it suffices to check bijectivity of $\theta_X$ just for free~$X$:
$X=\Sigma(S)=\Unit_\sA^{(S)}$. Since any such $X$ is a coproduct of a certain collection of copies of $\Unit_\sA$, and $F$ transforms coproducts into products, we have just to deal with the case $X=\Unit_\sA$. Then $\theta_{\Unit_\sA}:\Hom_\sA(\Unit_\sA,M)=M\to F(\Unit_\sA)=M$ is bijective (actually the identity map of $M$) by construction of~$\theta$.

This finishes the proof of completeness of $\catMod\Sigma$. In fact, we never used algebraicity or commutativity of monad~$\Sigma$ during this proof (these properties are needed for the proof of V1)--V5)), so we have actually proved the following:
\begin{Th} 
Category $\catMod\Sigma=\catSets^\Sigma$ of modules with respect to\/ {\em arbitrary\/} monad $\Sigma$ over $\catSets$ is complete.
\end{Th} 

\section{External tensor product of vectoids}

In this section we want to show how one can construct new vectoids from existing ones. The most interesting for us is the notion of {\em external tensor product\/} $\sA\boxtimes\sB$ of two vectoids $\sA$ and $\sB$, which, apart from other things, is the 2-product of $\sA$ and $\sB$ in the 2-category of vectoids. Such 2-products of vectoids are even easier to construct than 2-products of topoi.

\nxpoint (Analogy between vectoids and vector spaces.)
Before going on, we would like to present a certain informal point of view on vectoids, useful for understanding further discussions. Namely, one can think of a vectoid (without $\otimes$; cf.\ \ptref{p:vectoid.term}) as a counterpart, or a categorification, of a vector space. Cocontinuous functors between vectoids are then counterparts of linear maps of vector spaces, and bicocontinuous bifunctor $\otimes:\sA\times\sA\to\sA$ is an analogue of bilinear product on a vector space. Thus a $\otimes$-vectoid turns out to be an analogue of an algebra or a ring, and inverse image functors between $\otimes$-vectoids correspond to algebra homomorphisms. Furthermore, sets correspond to (small) categories in this picture, and the counterpart of a vector space with base $S$ is the presheaf vectoid $\hat\sS$, since any functor $\sS\to\sA$ from a small category into a vectoid extends essentially uniquely to a cocontinuous functor $\hat\sS\to\sA$.

\nxpoint (Origin of word ``vectoid''.)
The above analogy was actually the motivation for our terminology: vectoid is called so because it is somewhat similar to a vector space. However, in this work we use the word `vectoid' for the notion that would be more natural to call $\otimes$-vectoid (or {\em algebroid\/}, but this term already means something completely different). This somewhat weakens the original analogy (cf.\ \ptref{p:vectoid.term}).

\begin{Def}\label{def:tens.prod.vect} 
(External tensor product of vectoids.)
Let $\sA$, $\sB$ be two vectoids (not necessarily with $\otimes$).
Consider 2-category consisting of all bicocontinuous functors $F:\sA\times\sB\to\sC$, where $\sC$ is a variable vectoid (not necessarily with $\otimes$).
The initial object of this 2-category will be denoted by $\sA\times\sB\to\sA\boxtimes\sB$; vectoid (without $\otimes$) $\sA\boxtimes\sB$ will be called {\em the external tensor product of $\sA$ and~$\sB$}.
\end{Def}

Of course, this definition is the categorification of the definition of tensor product of vector spaces. Existence of $\sA\boxtimes\sB$ (uniqueness of which follows from 2-universal property) will be demonstrated only in the special cases we need further; nevertheless, the construction of $\sA\boxtimes\sB$ will be given for the general case as well.

\nxpoint (Case $\sC=\catSets^{op}$.)
Before proceeding with the construction of $\sA\boxtimes\sB$, we want to consider a special case. Namely, by definition the category of bicocontinuous bifunctors $\sA\times\sB\to\sC$ must be equivalent to the category of cocontinuous functors $\sA\boxtimes\sB\to\sC$ for any vectoid $\sC$ (without $\otimes$). Put here $\sC=\catSets^{op}$, regardless of the fact that this category isn't a vectoid (condition V4) fails). We obtain $\catFunct_{bicocont}(\sA\times\sB,\catSets^{op})\cong
\catFunct_{cocont}(\sA\boxtimes\sB,\catSets^{op})\cong\catFunct_{cont}((\sA\boxtimes\sB)^{op},\catSets)^{op}=\widetilde{\sA\boxtimes\sB}^{op}=(\sA\boxtimes\sB)^{op}$, since $\sA\boxtimes\sB$ is complete. Therefore, if the substitution of $\catSets^{op}$ for $\sC$ had been allowed, then $\sA\boxtimes\sB\cong\catFunct_{bicont}(\sA^{op}\times\sB^{op},\catSets)$. This suggests us that $\sA\boxtimes\sB$ might be constructed as the category of presheaves $F:\sA^{op}\times\sB^{op}\to\catSets$, continuous in each argument.

\nxpoint (Construction of $\sA\boxtimes\sB$.)
For any two vectoids (without $\otimes$) put
\begin{equation}
\sA\boxtimes\sB = \catFunct_{bicont}(\sA^{op}\times\sB^{op},\catSets)\subset
\widehat{\sA\times\sB}
\end{equation}
If $\sG_\sA\subset\sA$ and $\sG_\sB\subset\sB$ are small generating subcategories, then any bicontinuous functor $F:\sA^{op}\times\sB^{op}\to\catSets$ is completely determined by its restriction to $\sG_\sA\times\sG_\sB$; thus $\sA\boxtimes\sB$ may be identified with a full subcategory of presheaf category $\widehat{\sG_\sA\times\sG_\sB}$, hence $\sA\boxtimes\sB$ is well-defined from the set-theoretical point of view (e.g.\ it is a $\cU$-category with respect to chosen universe $\cU$).

\nxsubpoint (Alternative expressions for $\sA\boxtimes\sB$.)
Identify $\widehat{\sA\times\sB}=\catFunct(\sA^{op}\times\sB^{op},\catSets)$
with $\catFunct(\sA^{op},\catFunct(\sB^{op},\catSets))=\catFunct(\sA^{op},\hat\sB)=\hat\sA_{\hat\sB}$, i.e.\ with the category of $\hat\sB$-valued presheaves on $\sA$.
It is easy to see that $\sA\boxtimes\sB$ will be identified then with $\catFunct_{cont}(\sA^{op},\tilde\sB)=\catFunct_{cont}(\sA^{op},\sB)=\tilde\sA_\sB$, i.e.\ with the category of $\sB$-valued continuous presheaves on~$\sA$. Interchanging $\sA$ and $\sB$, we see that $\sA\boxtimes\sB\cong\tilde\sB_\sA$.

\nxsubpoint (Construction of bifunctor $\sA\times\sB\to\sA\boxtimes\sB$.)
For any continuous presheaves (``sheaves'') $X\in\Ob\tilde\sA$ and $Y\in\Ob\tilde\sB$ we define $X\boxtimes Y\in\Ob\widehat{\sA\times\sB}$ as follows:
\begin{equation}
X\boxtimes Y: (A,B)\mapsto X(A)\times Y(B)
\end{equation}
It is easy to see that $X\boxtimes Y$ is almost continuous in $A$ and in~$B$ (it preserves filtered colimits and cocartesian squares in each argument; the only obstruction to its bicontinuity is the initial object: $(X\boxtimes Y)(\emptyset,B)=Y(B)\neq\st1$). Denote by $X\tilde\boxtimes Y\in\Ob\sA\boxtimes\sB$ the associated bicontinuous functor (i.e.\ the value on $X\boxtimes Y$ of functor $a$, left adjoint to embedding $\sA\boxtimes\sB\to\widehat{\sA\times\sB}$); in this special case $X\tilde\boxtimes Y$ might admit an explicit description, because $X\boxtimes Y$ already has almost all necessary properties. This construction is clearly functorial in $X$ and~$Y$, hence it yields a functor $\Theta:\sA\times\sB=\tilde\sA\times\tilde\sB\to\sA\boxtimes\sB$.

\nxsubpoint (Complicated points in the proof.)
We will not verify $\sA\boxtimes\sB$ to be a vectoid, and $\Theta:\sA\times\sB\to\sA\boxtimes\sB$ to be the universal bicocontinuous functor, because the general construction of $\boxtimes$ won't be used in this work. We content ourselves by remarking that the most complicated points are the verification of V1) (existence of arbitrary colimits; here one will have to pass from a presheaf on $\sA\times\sB$ to associated bicocontinuous presheaf) and V6) (completeness, or totality) for $\sA\boxtimes\sB$, and then --- verification of bicocontinuity of bifunctor~$\Theta$.

\nxpoint (Special case: $\sA=\hat\sS$.)
When $\sA=\hat\sS$ is a presheaf category, any functor $\sS^{op}\to\sB$ can be uniquely extended (up to an isomorphism) to a cocontinuous functor $\sA^{op}\to\sB$ (since any presheaf is a colimit of some diagram of representable presheaves); hence $\hat\sS\boxtimes\sB=\catFunct_{cont}(\hat\sS^{op},\sB)=\catFunct(\sS^{op},\sB)=\hat\sS_\sB$ is the category of $\sB$-presheaves on $\sS$. In this special case it is much easier to check that $\sA\boxtimes\sB=\hat\sS_\sB$ is a vectoid and that $\sA\times\sB\to\sA\boxtimes\sB$ has the required universal property. We will use external tensor product of vectoids only in this case.

\nxsubpoint (Examples.)
In particular, categories of (co)simplicial objects in~$\sA$ are given by $s\sA=s\catSets\boxtimes\sA$ and $c\sA=c\catSets\boxtimes\sA$. Therefore, it would suffice to construct ``nice'' vectoids $s\catSets$ and $c\catSets$ once (for example, by convenient universal properties), and extend the construction to arbitrary vectoids by means of external tensor product afterwards.

\nxsubpoint (Bisimplicial sets.)
When $\sA=\hat\sS$ and $\sB=\hat\sT$ both are presheaf categories, then $\sA\boxtimes\sB=\widehat{\sS\times\sT}$. For example, $s\catSets\boxtimes s\catSets=s^2\catSets$ is the category of bisimplicial sets.

\nxpoint (2-$\otimes$-structure on 2-category of vectoids.)
Let us provisionally denote by $\catVectoid$ the 2-category of vectoids without tensor product, with cocontinuous functors as morphisms. Then the external tensor product 2-functor $\boxtimes:\catVectoid\times\catVectoid\to\catVectoid$ defines a 2-tensor structure on this category. Algebras (commutative associative with unit) in $\catVectoid$ are precisely the $\otimes$-vectoids (i.e.\ vectoids as defined in~\ptref{def:vectoid}), since giving a bicocontinuous functor $\otimes:\sA\times\sA\to\sA$ is equivalent to defining a morphism $\sA\boxtimes\sA\to\sA$ in $\catVectoid$. In other words, $\otimes$-vectoids are simply the commutative algebras inside $\catVectoid$ (however, with morphisms going in the opposite direction).

\nxsubpoint (Consequences for $\otimes$-vectoids.)
Similarly to the fashion in which tensor product of algebras turns out to be an algebra, $\sA\boxtimes\sB$ turns out to be a $\otimes$-vectoid any time both $\sA$ and $\sB$ are $\otimes$-vectoids themselves; and, similarly to the fact that tensor product of commutative algebras is their coproduct in the category of commutative algebras, $\sA\boxtimes\sB$ is the 2-product of $\sA$ and $\sB$ inside 2-category of $\otimes$-vectoids.

\section{Structure classifiers and Algebrads}

\nxpoint
Let $Z$ be a contravariant functor from 2-category of $\otimes$-vectoids into 2-category of categories. For example, $Z$ can transform a vectoid $\sA$ into the category of its objects (i.e.\ $Z(\sA)=\sA$), or into the category of associative algebras $\catAlg(\sA)$ in $\sA$ ($Z(\sA)=\catAlg(\sA)$, i.e.\ $Z=\catAlg$), or into the category of coalgebras $\catCoalg(\sA)$.

For any such $Z$ we can ask whether $Z$ is 2-representable by some $\otimes$-vectoid $\sC_Z$ together with some object $X_Z\in Z(\sC_Z)$. By definition, such representability means that for any $\otimes$-vectoid $\sA$ we have
\begin{equation}
Z(\sA)\cong\Hom(\sA,\sC_Z)=\catFunct_{\otimes,cocont}(\sC_Z,\sA)
\end{equation}
In other words, each object $X\in Z(\sA)$ appears as inverse image $f_X^*X_Z$ of the universal object $X_Z\in Z(\sC_Z)$ with respect to some uniquely determined (up to a unique isomorphism) $\otimes$-vectoid morphism $f_X:\sA\to\sC_Z$.

If such a vectoid $\sC_Z$ and object $X_Z\in Z(\sC_Z)$ do exist, we say that $\sC_Z$ is the {\em classifying vectoid} or the {\em classifier\/} of structures of type $Z$, and $X_Z$ is the {\em universal structure\/} of type~$Z$.

So, it makes sense to speak about the object classifier and the universal object, or about algebra classifier and universal algebra over it and so on.

\nxpoint
If we think of the classifying vectoid (say, of algebras) as a space, and of the universal structure (say, an algebra) on it as a certain sheaf over this space with some structure, then the points of this classifying vectoid (understood here as morphisms from the point vectoid $\catSets$ or, even better $\catVect k$, where $k$ is a field) correspond to elements $Z(\catVect k)$, i.e., say, to arbitrary algebras over arbitrary fields~$k$. In other words, the classifying vectoid can be understood as the {\em moduli space\/} for corresponding structures.

\nxpoint
While considering such problems it is sometimes convenient to work not only with ACU $\otimes$-vectoids, but with AU $\otimes$-vectoids as well, with not necessarily commutative tensor product. In such cases we will speak about ``not necessarily commutative'' classifying vectoids.

\nxpoint (Monads.)
Recall that a monad over some category $\sC$ can be defined as an algebra in the endofunctor category $\catEndof(\sC)$, considered as an AU $\otimes$-category with respect to the ``composition product'' $\otimes:=\circ$, given by composition of endofunctors.

One can replace $\sC$ with any object of any (strictly associative) 2-category. We will be interested in the case where $\sC$ is some $\otimes$-vectoid $\sZ$.

\begin{Def} (Definition of an algebrad.)
Let $\sZ$ be an arbitrary $\otimes$-vectoid. A {\em $\sZ$-algebrad}, or simply an {\em algebrad}, is by definition any algebra inside the AU $\otimes$-category of endomorphisms of vectoid $\sZ$ with respect to the composition product.
\end{Def}

We are especially interested in the case when $\sZ=\sC_Z$ is a classifying vectoid for sufficiently simple structures $Z$. Then the category $\catEnd_\catVectoid(\sC_Z)$, in which we are looking for algebras, admits a simpler description:
\begin{equation}\label{eq:endom.class.vect}
\catEnd_{\catVectoid}(\sC_Z)=\catHom_{\catVectoid}(\sC_Z,\sC_Z)=Z(\sC_Z)
\end{equation}
It is worthwhile to remark that $\catEnd_\catVectoid(\sC_Z)$ almost coincides with category $\catEndof_{\otimes,cocont}(\sC_Z)$: these two categories are equivalent, and differ only by the order of arguments of the composition product~$\circ$.

Equality \eqref{eq:endom.class.vect} often allows to obtain an explicit description of this category, describe its composition product $\circ$ and then describe algebras with respect to $\circ$, i.e.\ $\sC_Z$-algebrads. In this case we shall always denote the composition product by $\circ$, not by $\otimes$, in order not to confuse it with the tensor product of~$\sC_Z$.

\nxpoint\label{p:rel.comp.prod} ($\sZ$-algebrads over $\sA$.)
The definition of a $\sZ$-algebrad admits a relative version over any $\otimes$-vectoid $\sA$. Namely, consider $\sZ_\sA:=\sZ\boxtimes\sA$ as an object of the 2-category of vectoids over $\sA$. We still have $\catEnd_{\catVectoid_{/\sA}}(\sZ_\sA)=Z(\sZ_\sA)$, since $\sZ_\sA$ is a $Z$-structure classifier inside $\catVectoid_{/\sA}$. Now we can define a $\sZ$-algebrad $\Sigma$ over $\sA$ as an algebra in $\catEnd_\sA(\sZ_\sA)\cong Z(\sZ_\sA)$ with respect to composition product.

\nxpoint (Modules with respect to a $\sZ$-algebrad.)
Let $\Sigma$ be an arbitrary $\sZ$-algebrad over $\sA$, $\sB$ be any vectoid over $\sA$. In this situation we obtain a ``composition action''
\begin{equation}
\catEnd_\sA(\sZ_\sA)\times\catHom_\sA(\sB,\sZ_\sA)\to\catHom(\sB,\sZ_\sA)
\quad,
\end{equation} 
compatible with composition product in $\catEnd_\sA(\sZ_\sA)$, and since $\Sigma$ is an algebra in $\catEnd_\sA(\sZ_\sA)$ with respect to this product, it makes sense to speak about {\em $\Sigma$-modules\/} $M$ inside the category $\catHom(\sB,\sZ_\sA)$. For the sake of brevity we shall say that $M$ is a $\Sigma$-module over $\sB$, even if it actually lies in $\catHom(\sB,\sZ_\sA)$. If $\sZ=\sC_Z$ is a classifying vectoid, then $\catHom(\sB,\sZ_\sA)\cong Z(\sB)$, i.e.\ a $\Sigma$-module over $\sB$ actually lies in category $Z(\sB)$.

The definition of a $\Sigma$-module is most often applied for $\sB=\sA$.

\section{Object classifier and symmetric operads}

\zeropoint\nxpoint
The goal of this section is to construct the {\em classifying vectoid for objects\/}, or {\em object classifier\/}, which will be denoted by $\sQ_o$ or $q_o\catSets$. The universal object in $\Ob\sQ_o$ will be denoted by $X_o$.

\nxpoint (Notations.)
Let us fix following notations. For any integer $n\geq0$ denote by $\stn$ the standard finite set $\{1,2,\ldots,n\}$. Denote by $\catN$ the full subcategory of $\catSets$ consisting of all standard finite sets $\st0$, $\st1$, \dots, $\stn$, \dots; this category is equivalent to the category of finite sets $\catSets_f$.

Denote by $\catN_{iso}$ and $\catSets_{f,iso}$ subcategories of corresponding categories with the same objects, but with only bijective maps for morphisms. Thus $\Hom_{\catN_{iso}}(\stm,\stn)=\emptyset$ for $m\neq n$, and $\End_{\catN_{iso}}(\stn)=\gS_n$ is the symmetric group.

\nxpoint (Objects and morphisms obtained starting from an arbitrary object.)
Fix arbitrary $\otimes$-vectoid $\sA$ and its object $X\in\Ob\sA$. Which objects and morphisms in $\sA$ can be constructed starting only from~$X$? First of all, one can compute multiple tensor products of several copies of $X$, thus obtaining objects $X^{\otimes n}=X\otimes X\otimes\cdots\otimes X$ for all $n\geq0$. Furthermore, since $\otimes$ is an associative and commutative tensor structure on~$\sA$, the symmetric group $\gS_n$ acts on each $X^{\otimes n}$ by permuting factors. In this way we obtain a functor $\Phi_X:\catN_{iso}^{op}\to\sA$, transforming $\stn$ into $X^{\otimes n}$; this functor is completely determined by object $X$ up to an isomorphism.

\nxpoint\label{p:bar.phiX} (Extension to $\widehat{\catN_{iso}}$.)
Provisionally put $\sC:=\catN_{iso}$. Since $\sA$ admits arbitrary colimits, functor $\Phi_X:\sC\to\sA$, as well as any other functor, extends to a cocontinuous functor $\bar\Phi_X:\hat\sC\to\sA$, given by the usual formula
\begin{equation}\label{eq:cocont.ext}
\bar\Phi_X(F)=\bar\Phi_X(\injlim_{\sC_{/F}}p)=\injlim_{\sC_{/F}}\Phi_X\circ p,
\quad p:\sC_{/F}\to\sC
\end{equation}

\nxpoint ($\otimes$-structure on $\catN_{iso}$ and $\catSets_{f,iso}$.)
Moreover, $\sC$ admits a natural ACU $\otimes$-structure, given by disjoint union of finite sets: $\stn\otimes\stm:=\stn\sqcup\stm\cong \st{n+m}$. Since we want to work inside $\sC=\catN_{iso}$, we have to fix some bijection between this finite set and $\st{n+m}$, i.e\ some linear order on set $\stn\sqcup\stm$; we will usually assume that all elements of $\stn$ precede all elements of $\stm$ in this union. Under this approach the commutativity constraint $\stn\sqcup\stm\cong\stm\sqcup\stn$ corresponds to a non-trivial element of $\gS_{m+n}$.

Another option is to replace $\sC=\catN_{iso}$ by equivalent category $\catSets_{f,iso}$; then the disjoint union of finite sets can be computed in the natural fashion.

\nxpoint (Extending $\otimes$-structure to $\hat\sC$.)
Bifunctor $\otimes:\sC\times\sC\to\sC\subset\hat\sC$ can be uniquely extended to a bicocontinuous functor $\otimes:\hat\sC\times\hat\sC\to\hat\sC$, obviously defining an ACU $\otimes$-structure on $\hat\sC$. Thus $\hat\sC$ becomes a $\otimes$-vectoid. Besides, functor $\Phi_X:\sC\to\sA$, $\stn\mapsto X^{\otimes n}$ is compatible with tensor products of $\sC$ and $\sA$, since $X^{\otimes m}\otimes X^{\otimes n}\cong X^{\otimes(m+n)}$. By cocontinuity this compatibility extends to $\hat\sC$: functor $\bar\Phi_X:\hat\sC\to\sA$ is a $\otimes$-functor.

\nxpoint (Object classifier.)
Put $\sQ_o:=\hat\sC=\widehat{\catN_{iso}}$ together with the tensor product just defined, $X_o:=\st1$ (presheaf represented by standard one-element set). We have just shown that any object $X$ of any $\otimes$-vectoid $\sA$ defines some uniquely determined cocontinuous $\otimes$-functor $\bar\Phi_X:\sQ_o\to\sA$ transforming $X_o$ into $X$. This means exactly that $(\sQ_o,X_o)$ is the object classifier.

\nxpoint (Relative version.)
For any vectoid $\sA$ the object classifier inside  $\catVectoid_{/\sA}$ can be obtained by base change from $\sQ_o$: $\sQ_{o,\sA}=\sQ_o\boxtimes\sA=\widehat{\catN_{iso,\sA}}=\catFunct(\catN_{iso}^{op},\sA)$. Similarly to our notations for simplicial objects, put $q_o\sC:=\catFunct(\catN_{iso}^{op},\sC)$ for any category~$\sC$; then $\sQ_o=q_o\catSets$ and $\sQ_{o,\sA}=q_o\sA$.

\nxpoint (Explicit description of $q_o\sA$.)
Let us describe explicitly ``$q_o$-objects of category $\sA$'', i.e.\ the objects of $q_o\sA=\catFunct(\catN_{iso}^{op},\sA)$. An object $A\in\Ob q_o\sA$ is simply a sequence $(A_n)_{n\geq 0}$ of objects of $\sA$, together with right action of symmetric group $\gS_n$ on $A_n$ for each $n\geq0$. Morphism of $q_o$-objects $f:A\to B$ is a family of morphisms $f=(f_n)_{n\geq0}$, $f_n:A_n\to B_n$, compatible with all $\gS_n$-actions. 

This description is applicable in particular to the category of $q_o$-sets $q_o\catSets$: its objects are just collections of sets $A=(A_n)_{\geq0}$ with right $\gS_n$-action on~$A_n$.

\nxpoint (Representable objects in $\sQ_o=q_o\catSets$.)
Recall that representable objects $h_\stn=X_o^{\otimes n}$ constitute a full subcategory inside $\sQ_o$, equivalent to $\catN_{iso}$. In terms of $q_o$-sets this $h_\stn$ corresponds to sequence $(\emptyset,\ldots,\emptyset,\gS_n,\emptyset,\ldots)$. For example, $X_o$ corresponds to $(\emptyset,\st1,\emptyset,\ldots)$.

\nxpoint (Explicit description of $\bar\Phi_X:q_o\catSets\to\sA$.)
Let $X\in\Ob\sA$ be an arbitrary object of arbitrary $\otimes$-vectoid $\sA$. Let's describe cocontinuous $\otimes$- $\bar\Phi_X:q_o\catSets\to\sA$ of~\ptref{p:bar.phiX} explicitly. An easy computation using~\eqref{eq:cocont.ext} yields
\begin{equation}\label{eq:expl.phiX}
\bar\Phi_X(A)=\bigsqcup_{n\geq0}(A_n\otimes X^{\otimes n})/\gS_n
\end{equation}

\nxpoint (Tensor product of $q_o\catSets$ and $q_o\sA$.)
Using the same formulas together with equality $h_\stn\otimes h_\stm=h_{\st{m+n}}$ we obtain a formula for tensor product in $q_o\sA$:
\begin{equation}
(A\otimes B)_n=\bigsqcup_{p+q=n}\Ind_{\gS_p\times\gS_q}^{\gS_n}A_p\otimes B_q
\end{equation}
Here $\Ind_H^G(M)=M\otimes_H G$ for any groups $H\subset G$ and any object $M$ with right $H$-action. Formulas for multiple tensor products are similar:
\begin{equation}\label{eq:mult.tens.prod.qo}
(A^{(1)}\otimes\cdots\otimes A^{(s)})_n=\bigsqcup_{p_1+\cdots+p_s=n}\Ind^{\gS_n}_{\gS_{p_1}\times\cdots\times\gS_{p_s}}A^{(1)}_{p_1}\otimes\cdots\otimes A^{(s)}_{p_s}
\end{equation}

\nxpoint (Composition product on $q_o\sA\cong\catEnd_\sA(q_o\sA)$.)
Let us compute the {\em composition product\/} $\circ$ on category $\catEnd_\sA(q_o\sA)\cong q_o\sA$ (cf.\ \ptref{p:rel.comp.prod}).
Equivalence of categories $\theta:\catEnd_\sA(q_o\sA)\to q_o\sA$ transforms vectoid morphism $f=(f^*,f_*)$ with cocontinuous inverse image $\otimes$-endofunctor $f^*$ into its value $f^*(X_{o,\sA})$ on universal object $X_{o,\sA}=(\emptyset_\sA,\Unit_\sA,\emptyset_\sA,\ldots)$. Quasi-inverse equivalence $\theta^{-1}$ transforms $X$ into vectoid endomorphism with inverse image functor $\bar\Phi_X$, given by formula~\eqref{eq:expl.phiX} (or rather its relative version over $\sA$). One has also to bear in mind that the composition product is defined with respect to vectoid morphisms, not cocontinuous inverse image functors which are pointed in the opposite direction; hence $X\circ Y=\theta(\theta^{-1}(X)\circ\theta^{-1}(Y))=(\bar\Phi_Y\circ\bar\Phi_X)(X_{o,\sA})=\bar\Phi_Y(\bar\Phi_X(X_{o,\sA}))=\bar\Phi_Y(X)$, which can be computed using \eqref{eq:expl.phiX} and~\eqref{eq:mult.tens.prod.qo}:
\begin{equation}
X\circ Y=\bigsqcup_{s\geq0}(X_s\otimes Y^{\otimes s})/\gS_s\quad,
\end{equation}
whence
\begin{equation}\label{eq:compos.prod.qo}
(X\circ Y)_n=\bigl(\bigsqcup_{p_1+\cdots+p_s=n}X_s\otimes \Ind^{\gS_n}_{\gS_{p_1}\times\cdots\times\gS_{p_s}}(Y_{p_1}\otimes\cdots\otimes Y_{p_s})\bigr)/\gS_s\quad.
\end{equation}
Notice that the unit of $q_o\sA$ with respect to the tensor product is $\Unit_\otimes=(\Unit_\sA,\emptyset,\emptyset,\ldots)$, while the unit with respect to the composition product coincides with the (relative) universal object: $\Unit_\circ=\theta(\Id_{q_o\sA})=\Id_{q_o\sA}(X_{o,\sA})=X_{o,\sA}=(\emptyset,\Unit_\sA,\emptyset,\ldots)$. It's important to distinguish them.

\nxpoint ($\sQ_o$-algebrads over $\sA$.)
Now we are able to study $\sQ_o$-al\-geb\-rads over $\sA$, i.e.\ algebras $\Sigma=(\Sigma,\mu,\epsilon)$ with respect to composition product $\circ$ on $\catEnd_\sA(q_o\sA)\cong q_o\sA$. We have to define some object $\Sigma=(\Sigma_0,\Sigma_1,\ldots)$ from $q_o\sA$, i.e.\ sequence of objects $\Sigma_n$ of $\sA$ with right $\gS_n$-action. Next, giving the unit $\epsilon:\Unit_\circ=X_{o,\sA}\to\Sigma$ amounts to giving a morphism $\epsilon_1:\Unit_\sA\to\Sigma_1$. Defining the multiplication $\mu:\Sigma\circ\Sigma\to\Sigma$ is equivalent by~\eqref{eq:compos.prod.qo} to defining a family of $(\gS_{p_1}\times\cdots\times\gS_{p_s})$-equivariant morphisms in $\sA$:
\begin{equation}\label{eq:qo.alg.comb.descr}
\mu_{p_1,\ldots,p_s}:\Sigma_s\otimes\Sigma_{p_1}\otimes\cdots\otimes\Sigma_{p_s}\to\Sigma_{p_1+\cdots+p_s}
\end{equation}
Associativity and unit axioms can be also written in terms of these maps $\mu_{p_1,\ldots,p_s}$ and $\epsilon_1$. In the end we arrive, as the reader has already guessed, to the classical notion of a {\em symmetric operad over $\sA$}, which can be applied, for example, to $\sA=\catVect k$:

\begin{Th}
Let $\sQ_o$ be the classifying vectoid for objects. Then $\sQ_o$-algebrads over arbitrary vectoid $\sA$ are precisely the symmetric operads in $\sA$.
\end{Th}

\nxpoint
We want to emphasize that similar descriptions of symmetric operads (using an explicitly constructed category $q_o\sA$ and a composition product on it) are known for quite a long time. However, the classifying vectoids point of view helps us to understand where this category and composition product come from, and suggests to us possible further generalisations.

\nxpoint
Notice that $n$-sorted operads appear in our approach as algebrads, related to the classifier $\sQ_o\boxtimes\sQ_o\boxtimes\cdots\boxtimes\sQ_o=q_o^n\catSets$ of $n$-tuples of objects.

\nxpoint
As a final remark we would like to note that non-symmetric operads appear in a similar fashion, if one considers the object classifier inside the category of not necessarily commutative $\otimes$-vectoids.

\section{Coalgebra classifier and algebraic monads}

\zeropoint\nxpoint
Next we want to construct the {\em cocommutative coalgebra classifier $\sQ_c$} together with universal cocommutative coalgebra $X_c$ in $\sQ_c$, and study $\sQ_c$-algebrads. 

\nxpoint (Objects and morphisms constructed from a coalgebra.)
So, let us fix some $\otimes$-vectoid $\sA$ and cocommutative coalgebra $C=(C,\Delta,\eta)$, $\Delta:C\to C\otimes C$, $\eta:C\to\Unit$ in~$\sA$. Which objects and morphisms can be constructed in $\sA$ starting from this coalgebra?

\nxsubpoint (Morphisms corresponding to bijections in $\catN$.)
First of all, we still have tensor powers $C^{\otimes n}$ with right $\gS_n$-action on each; denote by $\sigma^*:C^{\otimes n}\simto C^{\otimes n}$ the automorphism induced by permutation $\sigma\in\gS_n=\Aut_\catN(\stn)$, then $(\sigma\circ\tau)^*=\tau^*\circ\sigma^*$. However, the counit and the comultiplication define many other morphisms; let's describe them.

\nxsubpoint\label{sp:coalg.monot.maps}
(Morphisms corresponding to non-decreasing maps in $\catN$.)
Denote by $\Delta^{(n)}:C\to C^{\otimes n}$ the {\em $n$-fold coproduct\/} morphism; thus $\Delta^{(0)}=\eta$, $\Delta^{(1)}=\id_C$, $\Delta^{(2)}=\Delta$, $\Delta^{(n+1)}=(\Delta^{(n)}\otimes\id_C)\circ\Delta$. Coassociativity and counit relations for $C$ yield $(\Delta^{(p)}\otimes\Delta^{(q)})\circ\Delta=\Delta^{(p+q)}$ for any $p$, $q\geq0$.

Now let $\phi:\stm\to\stn$ be an arbitrary non-decreasing map. Such a map is determined by partition of non-negative integer $m$ into $n$ non-negative summands $p_i:=\card\phi^{-1}(i)$; put
\begin{equation}
\phi^*:=\Delta^{(p_1)}\otimes\Delta^{(p_2)}\otimes\cdots\otimes\Delta^{(p_n)}:C^{\otimes n}\to C^{\otimes m}
\end{equation}
Coassociativity and counit relations imply that if $\psi:\stn\to\stp$ is another non-decreasing map, then $(\psi\circ\phi)^*=\phi^*\circ\psi^*:C^{\otimes p}\to C^{\otimes m}$. Moreover, cocommutativity of $C$ implies that for any permutation $\sigma:\stm\simto\stm$, such that $\phi=\phi\circ\sigma$, we have $\sigma^*\circ\phi^*=\phi^*$.

\nxsubpoint (Morphisms corresponding to arbitrary morphisms in $\catN^{op}$.)
Now let $f:\stm\to\stn$ be an arbitrary map of standard finite sets. It can be written as $f=\phi\circ\sigma$, where $\sigma\in\gS_m$ and $\phi:\stm\to\stn$ is a non-decreasing map. Put $f^*:=\sigma^*\circ\phi^*$. It is easy to see that this $f^*:C^{\otimes n}\to C^{\otimes m}$ doesn't depend on the choice of decomposition of~$f$; one also checks that $(g\circ f)^*=f^*\circ g^*$ for any $f:\stm\to\stn$ and $g:\stn\to\stp$. 

We have thus defined a functor $\Phi_C:\catN^{op}\to\sA$, transforming $\stn$ into $C^{\otimes n}$, and $f:\stm\to\stn$ into $f^*:C^{\otimes n}\to C^{\otimes m}$. Denote by $\bar\Phi_C:\hat\sC\to\sA$ its cocontinuous extension from $\sC$ to~$\hat\sC$, where $\sC:=\catN^{op}$.

\nxsubpoint (End of construction of $\sQ_c=\widehat{\catN^{op}}$.)
Notice that the disjoint union of sets defines on $\sC=\catN^{op}$ a tensor structure, with respect to which $\Phi_C:\sC\to\sA$ becomes a $\otimes$-functor. Define $\otimes$-structure on $\sQ_c:=\hat\sC$ by means of bicocontinuous extension of $\otimes:\sC\times\sC\to\sC$; then $(\sQ_c,\otimes)$ becomes a $\otimes$-vectoid, and $\bar\Phi_C:\hat\sC\to\sA$ a cocontinuous $\otimes$-functor, i.e.\ an inverse image functor for some vectoid morphism.

It remains to note that by construction $\bar\Phi_C$ transforms representable object $X_c:=h_{\st1}$ into $C$; morphisms $\Delta_c:=u^*$ and $\eta_c:=v^*$, where $u:\st2\to\st1$ and $v:\st0\to\st1$ are the only possible maps, are transformed into the comultiplication and counit of $C$, and themselves define a cocommutative coalgebra $X_c=(X_c,\Delta_c,\eta_c)$ inside $\sQ_c$. This proves that $\sQ_c$ is the cocommutative coalgebra classifier, and $X_c$ is the universal cocommutative coalgebra.

\nxpoint (Description of $\sQ_c$.)
Category $\sQ_c=\catFunct(\catN,\catSets)$ and its tensor product admit several different descriptions. First of all, notice that representable functor $h_\stn=\Hom_\catSets(\stn,-)$ transforms any finite set $S$ into $S^n$; this immediately implies that $h_\stn\otimes h_\stm=h_{\st{m+n}}$ is isomorphic to direct product $h_\stn\times h_\stm$. Since the direct product on $\sQ_c$ is bicocontinuous and coincides with the tensor product on $\catN^{op}$, they have to coincide everywhere: $\otimes_{\sQ_c}=\times$. In other words, {\em vectoid $\sQ_c$ is a topos, namely, the topos of presheaves on $\catN^{op}$}.

\nxsubpoint (Relation to algebraic endofunctors.)
Let $\catEndof_{alg}(\catSets)$ be the full subcategory of $\catEndof(\catSets)$ consisting of those endofunctors $F:\catSets\to\catSets$, which preserve filtered inductive limits. Since any set is a filtered inductive limit of finite sets, and category $\catSets_f$ is equivalent to $\catN$, it is easy to see that restriction to $\catN$ determines an equivalence of categories $\catEndof_{alg}(\catSets)\cong\catFunct(\catN,\catSets)$; quasi-inverse equivalence is given then by the left Kan extension along $\catN\to\catSets$.

In other words, {\em category $\sQ_c$ is equivalent to the category of algebraic endofunctors on $\catSets$.} Since representable functors $h_\stn$ are algebraic, it is easy to see that under this equivalence tensor product $\otimes_{\sQ_c}$ still corresponds to direct product of algebraic endofunctors.

\nxsubpoint (Description of $\bar\Phi_X(F)$ for any $F$.)
Let $F$ be an arbitrary object of $\sQ_c$, identified with corresponding algebraic endofunctor $F:\catSets\to\catSets$, $X$ be any cocommutative coalgebra in $\sA=\catSets$, i.e.\ any set (endowed with its only coalgebra structure with respect to the direct product of $\catSets$). Our general considerations imply that $\bar\Phi_X(F)$ may be computed as follows:
\begin{equation}\label{eq:expl.phiX.qc}
\bar\Phi_X(F)=\injlim_{\catN^{op}_{/F}}X^n=\injlim_{\stn\to F}X^n=
\Coker\bigl(\bigsqcup_{\phi:\stm\to\stn}F(m)\times X^n\rightrightarrows\bigsqcup_{n\geq0}
F(n)\times X^n\bigr)
\end{equation}
On the other hand, this same cokernel coincides with colimit $\injlim_{\stn\to X}F(n)$, which computes the value $F(X)$ of algebraic endofunctor $F$ on any set~$X$. Therefore,
\begin{equation}
\bar\Phi_X(F)=F(X)
\end{equation}

\nxpoint (Description of composition product in $\sQ_c$.)
The formula just proved enables us to compute the composition product in $\sQ_c$: it turns out to coincide with composition of algebraic endofunctors (taken in opposite order, which is not so important now). Therefore,
$\sQ_c$-algebrads over $\catSets$ are identified with algebraic monads over $\catSets$ in the sense of~\cite[4]{Du}, or, which is essentially the same, with algebraic theories (cf.~\cite{Lawvere}). Modules in $\catSets$ with respect to a $\sQ_c$-algebrad are precisely modules over an algebraic monad, or models of an algebraic theory.

\nxpoint (Combinatorial description.)
Of course, the composition product $F\circ G$ admits a combinatorial description in terms of sets $F(n)$ and $G(n)$, similar to~\eqref{eq:compos.prod.qo}:
\begin{equation}
(F\circ G)(n)=\Coker\Bigl(\bigsqcup_{\phi:\stp\to\stq}F(p)\times G(n)^q
\rightrightarrows\bigsqcup_{p\geq0}F(p)\times G(n)^p\Bigr)
\end{equation}
Therefore, the multiplication $\mu:\Sigma\circ\Sigma\to\Sigma$ of algebraic monad (i.e.\ $\sQ_c$-algebrad) $\Sigma$ may be described by a compatible system of maps of sets $\mu_{p,n}:\Sigma(p)\times\Sigma(n)^p\to\Sigma(n)$; this leads to the completely combinatorial description of algebraic monads given in~\cite[4]{Du}.

\nxsubpoint (Relative case.)
Notice that if $\Sigma$ is an arbitrary algebraic monad and $\sA$ is a vectoid, then a $\Sigma$-module over $\sA$ is not just an object of $\sA$, but a cocommutative coalgebra in $\sA$, endowed with $\Sigma$-action. When we work with topoi, and not just with arbitrary vectoids, we fail to see this additional coalgebra structure since any object of a topos admits a unique coalgebra structure with respect to $\otimes=\times$, given by the diagonal map. If we want to be able to define $\Sigma$-module structures on arbitrary objects of arbitrary vectoids, we have to consider operads instead of algebraic monads.

\nxpoint ($\sQ_c$-algebrads over $\sA$.)
Of course, it makes sense to consider $\sQ_c$-algebrads $\Sigma$ not only over $\catSets$, but over arbitrary vectoids $\sA$ as well. We have to use the combinatorial description in this case: $\Sigma$ will be a cocommutative coalgebra in $q_c\sA=\sQ_c\boxtimes\sA=\catFunct(\catN,\sA)$, i.e.\ apart from objects $\Sigma(n)\in\Ob\sA$ and morphisms $\phi_*:\Sigma(m)\to\Sigma(n)$ for each $\phi:\stm\to\stn$ we must require each $\Sigma(n)$ to be a cocommutative coalgebra in $\sA$, and all $\phi_*$ to be coalgebra homomorphisms (since the tensor product in $q_c\sA$ is computed componentwise: $(F\otimes G)(n)=F(n)\otimes G(n)$). After this $\epsilon:\Unit_\circ\to\Sigma$ and $\mu:\Sigma\circ\Sigma\to\Sigma$ admit a combinatorial description in terms of morphisms $\epsilon_1:\Unit_\sA\to\Sigma(1)$ and $\mu_{p,n}:\Sigma(p)\otimes\Sigma(n)^{\otimes p}\to\Sigma(n)$. Additional conditions must be imposed on these morphisms, completely similar to those arising over $\catSets$.

\nxpoint (Non-commutative case.)
In a similar, but easier, fashion one can construct the coalgebra classifier in the 2-category of not necessarily commutative (i.e.\ AU) $\otimes$-vectoids. In this case morphisms $\phi^*:C^{\otimes n}\to C^{\otimes m}$ will be defined for non-decreasing $\phi$ only, and we will finally obtain $\hat\Delta_\epsilon=s_\epsilon\catSets$, the category of presheaves over the category of finite ordinals (not necessarily non-empty), i.e.\ to the category of {\em augmented simplicial sets.} However, the tensor product doesn't coincide with the direct product in this case: since it is given by disjoint union of ordinals on simplices ($\Delta_\epsilon(n)\otimes\Delta_\epsilon(m)=\Delta_\epsilon(m+n+1)$), we see that $\otimes$ coincides with the {\em join\/} $\star$ of augmented simplicial sets.

\section{Algebra classifier}

\zeropoint\nxpoint
The goal of this section is to construct the {\em commutative algebra classifier\/ $\sQ_a$}. If $\sQ_o$ led us to symmetric operads, and $\sQ_c$ --- to algebraic monads, where will $\sQ_a$ lead us?

\nxpoint (Objects and morphisms constructed starting from an algebra.)
Let $\sA$ be an arbitrary vectoid again, $A=(A,\mu,\epsilon)$ be a commutative algebra in $\sA$. Construct objects $A^{\otimes n}$ as before, but this time we want $\gS_n$ to act from the {\em left\/} by means of automorphisms $\sigma_*:A^{\otimes n}\to A^{\otimes n}$; in the notations of previous section $\sigma_*=(\sigma^*)^{-1}$. Starting from $p$-fold product morphisms $\mu^{(p)}:A^{\otimes p}\to A$ we construct $\phi_*:A^{\otimes m}\to A^{\otimes n}$ for any non-decreasing map $\phi:\stm\to\stn$, using construction dual to~\ptref{sp:coalg.monot.maps}. Combining together $\phi_*$ and $\sigma_*$, define $f_*:A^{\otimes m}\to A^{\otimes n}$; in this way we obtain a functor $\Phi_A:\catN\to\sA$, which we extend to a cocontinuous functor $\bar\Phi_A:\hat\catN\to\sA$ in our usual fashion.

\nxpoint (End of construction of $\sQ_a:=\hat\catN$.)
Define $\otimes$ on $\catN$ by means of the disjoint union of finite sets, and extend this tensor product to all of $\sQ_a=\hat\catN$ by cocontinuity. Then $\sQ_a$ endowed with this tensor product turns out to be a vectoid, and $\bar\Phi_A$ turns out to be a inverse image functor for some vectoid morphism. Besides, the universal commutative algebra in~$\sQ_a$ is given by $X_a=(h_{\st1},\mu_a:=u_*,\epsilon_a:=v_*)$, $u:\st2\to\st1$, $v:\st0\to\st1$, and $\bar\Phi_A$ transforms $X_a$ into original commutative algebra $A$. This shows that $(\sQ_a,X_a)$ is the commutative algebra classifier.

\nxpoint (Explicit formula for $\bar\Phi_A(F)$.)
Let us present an explicit formula for $\bar\Phi_A(F)$, similar to \eqref{eq:expl.phiX} and~\eqref{eq:expl.phiX.qc}:
\begin{equation}\label{eq:expl.phiX.qa}
\bar\Phi_A(F)=\Coker\bigl(\bigsqcup_{\phi:\stm\to\stn}F(n)\otimes X^{\otimes m}\rightrightarrows\bigsqcup_{n\geq0} F(n)\otimes X^{\otimes n}\bigr)
\end{equation}

\nxpoint (Explicit formula for tensor product in $\sQ_a$.)
It is convenient to replace $\sQ_a=\hat\catN$ with equivalent category $\widehat{\catSets_f}$ for writing the next formula, i.e.\ we assume that all presheaves $F:\catN^{op}\to\catSets$ are naturally extended to $\catSets_f$. Then
\begin{equation}\label{eq:expl.otimes.qa}
(F\otimes G)(W)=\bigsqcup_{U\sqcup V=W}F(U)\otimes G(V),\quad
\text{where $U$, $V$, $W\in\Ob\catSets_f$}
\end{equation}
This formula is valid in $q_a\sA=\sQ_a\boxtimes\sA=\catFunct(\catN^{op},\sA)\cong\catFunct(\catSets_f^{op},\sA)$; of course, if we work in $\sQ_a=q_a\catSets$, the tensor product in the right-hand side becomes the direct product of sets.

In order to verify this formula it is sufficient to notice that its right-hand side defines a certain bicocontinuous tensor product on $q_a\sA$, coinciding with the required one on generators $X_a^{\otimes n}$. For example, for $\sA=\catSets$, $F=h_\stn$, $G=h_\stm$ we have $(h_\stn\otimes h_\stm)(W)=\Hom(W,\stn\sqcup\stm)$ on the left and $\sqcup_{U\sqcup V=W}\Hom(U,\stn)\times\Hom(V,\stm)$ on the right; these two expressions are easily seen to be canonically isomorphic.

Notice that $\Unit_\otimes=h_{\st0}$ is given by $\Unit_\otimes(0)=\Unit_\sA$, $\Unit_\otimes(n)=\emptyset$ for $n>0$ in this case.

\nxpoint\label{p:alg.qa} (Algebras in $\sQ_a$ and $q_a\sA$.)
In order to describe the category $\catEnd_\sA(q_a\sA)\cong\catCommAlg(q_a\sA)$, which contains the underlying objects for $\sQ_a$-algebrads, we have to study commutative algebras in $q_a\sA$. Let $(\Sigma,M:\Sigma\otimes\Sigma\to\Sigma,E:\Unit_\otimes\to\Sigma)$ be an algebra in $q_a\sA$. According to~\eqref{eq:expl.otimes.qa}, $M$ is given by a family of morphisms $M_{U,V}:\Sigma(U)\otimes\Sigma(V)\to\Sigma(U\sqcup V)$ for all finite sets $U$ and $V$; they can be also written $M_{p,q}:\Sigma(p)\otimes\Sigma(q)\to\Sigma(p+q)$ if one so wishes, but commutativity and associativity conditions for $M$ will become less natural. Unit $E$ is given by $E_0:\Unit_\sA\to\Sigma(\emptyset)$, as usual.

\nxpoint (Explicit formula for combinatorial product.)
Combining formulas \eqref{eq:expl.phiX.qa} and~\eqref{eq:expl.otimes.qa}, we can obtain an explicit formula for the composition product $X\circ Y=\bar\Phi_Y(X)$ of two commutative algebras in $q_a\sA$ (in fact, only $Y$ has to be a commutative algebra here):
\begin{gather}
\begin{split}
(X\circ Y)(V)=\Coker\Bigl(&\bigsqcup_{\phi:\stm\to\stn}\bigsqcup_{U_1\sqcup\cdots\sqcup U_m=V}X(n)\otimes Y(U_1)\otimes\cdots\otimes Y(U_m)\\
\rightrightarrows&
\bigsqcup_{n\geq0}\bigsqcup_{U_1\sqcup\cdots\sqcup U_n=V}X(n)\otimes Y(U_1)\otimes\cdots\otimes Y(U_n)\Bigr)
\end{split}
\end{gather}

\nxpoint (Combinatorial description of $\sQ_a$-algebrads over $\sA$.)
The above formulas enable us to obtain a combinatorial description of a $\sQ_a$-algebrad $\Sigma=(\Sigma,\mu,\epsilon)$ over $\sA$. First of all, $\Sigma$ lies in $\catEnd_\sA(q_a\sA)\cong\catCommAlg(q_a\sA)$. This means that some $\Sigma(U)\in\Ob\sA$ are given, depending contravariantly on finite set $U$, together with a compatible system of maps $M_{U,V}:\Sigma(U)\otimes\Sigma(V)\to\Sigma(U\sqcup V)$ and $E_0:\Unit_\sA\to\Sigma(\emptyset)$. Next, $\epsilon:\Unit_\circ=X_a\to\Sigma$ is given by $\epsilon_1:\Unit_\sA\to\Sigma(1)$, and the multiplication $\mu:\Sigma\circ\Sigma\to\Sigma$ can be described by means of maps
\begin{equation}
\mu_{p_1,\ldots,p_n}:\Sigma(n)\otimes\Sigma(p_1)\otimes\cdots\otimes\Sigma(p_n)\to
\Sigma(p_1+\cdots+p_n)
\end{equation}
These maps are similar to those arising in operad theory, i.e.\ for $\sQ_o$-algebrads (cf.\ \eqref{eq:qo.alg.comb.descr}), however, they are required to satisfy stronger compatibility relations. We will not write down explicitly all these relations between different $M_{p,q}$, $\mu_{p_1,\ldots,p_n}$, $E_0$ and $\epsilon_1$; however, all necessary tools for this are ready.

\nxpoint
Thus the classifying vectoid of commutative algebras defines a completely new notion of a $\sQ_a$-algebrad, which is somewhat similar to those of symmetric operads and algebraic theories, but is not reduced to them. Since these two other cases ($\sQ_o$ and $\sQ_c$) have led us to well-studied and definitely profound notions, it is natural to suppose that $\sQ_a$-algebrads should be very interesting as well, and that they merit a detailed study by themselves. We finish our brief excursion into the theory of classifying vectoids and corresponding algebrads at this point.

\bigbreak

\end{document}